\documentclass[number,seceqn,citesort,dvips]{arxbj}
\usepackage{cursive}  
\usepackage[vtexbibl]{}

\aid{0}
\volume{16}
\issue{4}
\pubyear{2010}
\firstpage{995}
\lastpage{1015}
\doi{10.3150/10-BEJ253}

\makeatletter

\newtheorem{lemma}{Lemma}[section]
\newtheorem{theorem}[lemma]{Theorem}
\newtheorem{corollary}[lemma]{Corollary}
\newremark{example}[lemma]{Example}
\newremark{remark}[lemma]{Remark}

\def\eqref#1{(\ref{#1})}
\makeatother

\begin{document}
\begin{frontmatter}

\title{Weak convergence of the function-indexed integrated periodogram for infinite variance processes}
\runtitle{Weighted periodogram for infinite variance processes}

\begin{aug}
\author[1]{\fnms{Sami Umut} \snm{Can}\thanksref{1}\ead[label=e1]{suc2@cornell.edu}},
\author[2]{\fnms{Thomas} \snm{Mikosch}\thanksref{2}\ead[label=e2]{mikosch@math.ku.dk}\corref{}} \and\break
\author[3]{\fnms{Gennady}~\snm{Samorodnitsky}\thanksref{3}\ead[label=e3]{gennady@orie.cornell.edu}}
\runauthor{S.U. Can, T. Mikosch and G. Samorodnitsky}
\address[1]{Center for Applied Mathematics,
Cornell University,
657 Rhodes Hall,
Ithaca, NY 14853, USA.\\ \printead{e1}}
\address[2]{Department of Mathematics,
University of Copenhagen,
Universitetsparken 5,
DK-2100 Copenhagen, Denmark. \printead{e2}}
\address[3]{School of Operations Research and Information Engineering,
Cornell University,
220 Rhodes Hall,
Ithaca, NY 14853, USA. \printead{e3}}
\end{aug}

\received{\smonth{3} \syear{2009}}
\revised{\smonth{9} \syear{2009}}

\begin{abstract}
In this paper, we study the weak convergence of the integrated periodogram
indexed by classes
of functions for linear processes with symmetric $\alpha$-stable
innovations. Under suitable summability conditions on the series of the
Fourier coefficients of the index functions, we show that the weak limits
constitute $\alpha$-stable processes which have representations as infinite
Fourier series with i.i.d.~$\alpha$-stable coefficients. The
cases $\alpha\in(0,1)$ and $\alpha\in[1,2)$ are dealt with by rather
different methods and under different assumptions on the classes of
functions. For example, in contrast to the case $\alpha\in(0,1)$,
entropy conditions are needed for $\alpha\in[1,2)$ to ensure the
tightness of the sequence of integrated periodograms indexed by
functions.
The results of this paper are of additional interest since
they provide limit results for infinite mean random quadratic forms
with
particular Toeplitz coefficient matrices.
\end{abstract}

%

\end{frontmatter}

\section{Introduction}\label{s0}

Over the last decades, efforts have been made to get a better understanding
of non-Gaussian time series in the time and frequency
domains.
In particular, time series whose marginal
distributions exhibit power law behavior have attracted a lot of attention.
The need for such models arises from
applications in areas as diverse as insurance, geophysics, finance and
telecommunications. Infinite fourth
moments are not untypical for series of daily log-returns from
exchange rates, stock indices, and other speculative prices, whereas
infinite second moments can be observed in time series from insurance
such as for windstorm, industrial fire and earthquake
insurance \cite{davismikosch2009a,davismikosch2009b,embrechts1997,mikosch2003}.
Infinite first moments are
typical for the marginal distribution of the magnitudes of
earthquakes \cite{kagan1997}. Infinite
variances are observed for the sizes of teletraffic data in
the World Wide Web \cite{crovellabestavros1996,crovellabestavrostaqqu1996,faygonzale2006,leland1993,willingeraqq1995};
see also the recent books \cite{adlerepsteintaqqu1998,resnick2006}.

Classical time series analysis deals with the second
(or higher) order moment structure of a stationary sequence. Heavy-tailed
modeling
requires, in addition, that one takes into account the interplay
between the
dependence structure and the tails of the series.
An important task is to understand the
classical statistical estimators and test procedures
when big shocks to the underlying system are present.
When the marginal distributions have infinite variance,
the notions of autocovariance, autocorrelation and spectral
distribution lose their meaning.
However, various studies over the last twenty years have shown
that the analysis of linear processes
$X_t=\sum_{j=-\infty}^{\infty} \psi_j  \varepsilon_{t-j}$,
$t\in{\mathbb{Z}}$,
with heavy-tailed i.i.d.~{\em innovations} $(\varepsilon_j)_{j\in
\mathbb{Z}}$
and constant coefficients $(\psi_j)_{j\in\mathbb{Z}}$ is very similar
to classical (finite variance) time series
analysis, where notions such as
autocovariances and spectral density are defined only in terms of
the $\psi_j$'s and the innovation variance $\sigma_\varepsilon^2$.
Most estimators and test statistics from classical time
series analysis can be modified insofar that one considers
self-normalized (or studentized) versions of them and for these
versions, an asymptotic
theory exists which parallels the classical theory with Gaussian
limit processes. In contrast to the latter theory,
the limits involve infinite variance stable distributions and
processes \cite{embrechts1997,kluppelbergmikosch1996b,mikosch1998}.

One of the main goals of classical (finite variance) time series
analysis is
to study the spectral properties of the linear process $(X_t)$.
In this context, the {\em periodogram}
\[
I_{n,X}(\lambda)=
\Biggl|\frac{1}{\sqrt{n}}\sum_{t=1}^n \mathrm{e}^{-\mathrm{i}\lambda t}
X_t\Biggr|^2 ,
\qquad\lambda\in[0,\curpi] ,
\]
plays a prominent role as an estimator
of the spectral density. Numerous estimation and test procedures are
based on this statistic and integrated versions of the form
$
J_{n,X}(f)=\int_0^{\curpi} I_{n,X}(\lambda) f(\lambda)\, \mathrm{d}\lambda
$
for appropriate
classes of real-valued functions $f\in\mathcal{F}$ on $[0,\curpi]$. In
applications, the class
$\mathcal{F}_I= \{I_{[0,x]}\dvt x\in[0,\curpi]\}$ is most important.
The resulting integrated periodogram is a process indexed by $x\in
[0,\curpi]$.
Under general conditions, $J_{n,X}(I_{[0,\cdot]})$ converges uniformly
with probability 1 to the function
$\sigma_\varepsilon^2 \int_{0}^\cdot|\psi(\mathrm{e}^{-\mathrm{i}\lambda
})|^2
\,\mathrm{d}\lambda$,
where
\begin{eqnarray*}
\psi(\mathrm{e}^{-\mathrm{i}\lambda})=\sum_{j=-\infty}^{\infty} \psi_j
\mathrm{e}
^{-\mathrm{i}\lambda
j} ,\qquad\lambda\in[0,\curpi] ,
\end{eqnarray*}
and
$|\psi(\mathrm{e}^{-\mathrm{i}\lambda})|^2$ is the corresponding {\em power transfer
function}. The latter is one of the essential building blocks of the
{\em spectral density} of the stationary process $(X_t)$:
\[
f_X(\lambda)=\dfrac{\sigma_\varepsilon^2}{2\curpi}
|\psi(\mathrm{e}^{-\mathrm{i}\lambda})|^2 =\frac{1}{2\curpi}\sum_{h=-\infty
}^{\infty}
\mathrm{e}^{-\mathrm{i}h\lambda}  \gamma_X(h) ,\qquad\lambda\in[0,\curpi] .
\]
This is the
Fourier series based on the {\em autocovariance function}
\[
\gamma_X(h)=\operatorname{cov}(X_0,X_h)=\sigma_\varepsilon^2 \sum
_{j=-\infty
}^{\infty}
\psi_j\psi_{j+|h|} ,\qquad h\in\mathbb{Z} .
\]
%
Since
$J_{n,X}(I_{[0,\cdot]})$ estimates the spectral distribution function of
the stationary process $(X_t)$, it has been used for a long time as the
{\em
empirical spectral distribution function}, both as an estimator and
a basic tool for constructing goodness-of-fit tests for the
underlying spectral distribution function
\cite{grenanderrosenblatt1984,brockwelldavis1991,priestley1981}.
Since the limit process of the properly centered and normalized
process $J_{n,X}(I_{[0,\cdot]})$ depends on the (in general unknown) spectral
density $f_X$, Bartlett \cite{bartlett1954} proposed to consider the
integrated periodogram\ based on
$\mathcal{F}_{\mathrm{B}}=\{I_{[0,x]}/f_X(x)\dvt x\in[0,\curpi]\}$.
Under general conditions, this process converges
uniformly, with probability 1, to the function $f(x)\equiv x$ and the
limit process can be shown to be independent of the coefficients of
the linear process,
but depends on the fourth moment of $\varepsilon_1$. More generally,
weighted integrated periodograms based on the classes
${\mathcal{F}}_g=\{I_{[0,x]} g(x)\dvt x \in[0,\curpi]\}$ for suitable
functions $g$
are used to estimate the spectral density or to perform
various tests concerning the spectrum of the underlying stationary
sequence \cite{priestley1981}.

The weighted integrated periodogram is also the basis of one of
the classical estimators for fitting ARMA and fractional ARIMA
models. This method goes back to early work by Whittle
\cite{whittle1951}. In this context, one considers the functional
$J_{n,X}(1/f_X(\cdot;\theta))$, $f_X(\cdot;\theta)\in
\mathcal{F}_\mathrm{W}$, where $\mathcal{F}_\mathrm{W}$ is a class of spectral
densities indexed by a parameter $\theta\in\Theta\subset
\mathbb{R}^d$. The {\em Whittle estimator} $\widehat\theta_n$ of the
true parameter $\theta_0 \subset\Theta$ is the minimizer of
$J_{n,X}(1/f_X(\cdot;\theta))$ over $\Theta$.
This estimation technique is
one of the backbones of quasi-maximum likelihood estimation in parametric
time series modeling. The Whittle estimator is known to be
asymptotically equivalent to the corresponding least-squares and
Gaussian quasi-maximum likelihood estimators \cite{brockwelldavis1991}.
When proving the
asymptotic normality and consistency of $\widehat\theta_n$,
one has to study the properties of the sequence
$(J_{n,X}(1/f_X(\cdot;\widehat\theta_n)))$ which, again, can be
considered as a weighted integrated periodogram indexed by a class of functions.

The above examples have in common
that one always considers a stochastic process $(J_{n,X}(f))_{f\in
\mathcal{F}}$ for some class of functions. In all cases, one is
interested in the asymptotic behavior of the process $J_{n,X}$,
uniformly over the class $\mathcal{F}$. This is analogous to the case
of the empirical distribution function indexed by classes of
functions. General references in this context are the monographs
\cite{pollard1984,vandervaartwellner1996}.
Early on, this analogy was discovered by Dahlhaus \cite{dahlhaus1988}, who gave
the uniform convergence theory for $J_{n,X}$ under entropy and
exponential moment
conditions. The almost sure and weak convergence
theory under entropy and power moment conditions
was given in \cite{mikoschnorvaisa1997}.
A recent survey of nonparametric statistical methods related to the
empirical spectral distribution
indexed by classes of functions is \cite{dahlhauspolonik2002}.

It is the aim of this paper to develop an analogous weak convergence
theory for heavy-tailed stationary processes. We will
understand `heavy-tailedness' in the sense of infinite variance
of the marginal distributions. Our focus will be on linear
processes $(X_t)$ with i.i.d.~symmetric
$\alpha$-stable (S$\alpha$S) innovations $(\varepsilon_t)$ for some
$\alpha\in(0,2)$.
Recall that a random variabl $Y_\alpha$ has a {\em symmetric stable
distribution}
$(Y_{\alpha} \sim S_{\alpha} (\sigma, 0, 0))$ if there are parameters
$0 < \alpha\leq2$, and $\sigma\geq0$ such that its characteristic
function has the form $E \mathrm{e}^{\mathrm{i} t Y_{\alpha}} = \mathrm{e}^{-\sigma
^\alpha
|t|^\alpha}$. For convenience, we also assume that \mbox{$\sigma=1$} for the
distribution
of $\varepsilon_1$. For $\alpha<2$, the random variable
$\varepsilon_1$ is known to have
infinite
variance \cite{feller1971,samorodnitskytaqqu1994}.
Much
of the theory given below depends on
tail estimates for random quadratic forms in i.i.d.~infinite variance
random variables. Such results are available for i.i.d.~stable
sequences. Although it seems
feasible that the theory can be extended to the more general class of
processes whose innovations have regularly varying\ tails, we do not
attempt to achieve this goal. The price would be more technicalities and
the gain would be negligible.

We intend to
show how the
classical (finite variance) tools and methods have to be modified in the
infinite variance stable situation
which can be considered as a boundary case of
the classical one when some of the innovations assume extremely
large values. By now, there exists quite a clear
picture concerning the asymptotic theory of the sample autocovariances, the
periodogram and its integrated versions when the innovation sequence\
in a
linear process has infinite variance; see
\cite{embrechts1997}, Chapter 7. In addition to
the latter reference,
goodness-of-fit tests for heavy-tailed processes
(corresponding to the class $\mathcal{F}_I$) were considered for
short- and long-memory linear processes
\cite{kluppelbergmikosch1996a,kokoszkamikosch1997}, and
Whittle estimation for infinite variance ARMA and FARIMA processes
was also studied
\cite{mikoschgadrich1995,kokoszkataqqu1996}.

The paper is organized as follows.
In Section \ref{sec:1}, we introduce
some useful notation for the integrated
periodogram. Our main goal is to prove the weak convergence of the integrated
periodograms indexed by suitable classes of functions. We
achieve this goal for an i.i.d.~sequence  in Section~\ref{sec:iid}, first
by showing the {\em convergence
of the finite-dimensional distributions} (Section~\ref{s1}), then
the {\em tightness}. The conditions and methods are rather different
in the cases $\alpha\in(0,1)$ (Section~\ref{s2})
and $\alpha\in[1,2)$ (Section~\ref{s3}).
The case $\alpha\in(0,1)$ is treated in the more general context of
random quadratic forms with Toeplitz coefficient matrices
satisfying some summability condition. The case $\alpha\in[1,2)$
requires entropy conditions and the corresponding techniques.
In Section~\ref{sextra}, we
extend the limit theory for the integrated periodograms from an
i.i.d.~sequence to linear processes.
The \hyperref[app]{Appendix} contains some auxiliary results concerning tail estimates of
random quadratic forms in stable random variables. The weak
convergence results of this
paper might also be of
separate interest in the context of infinite variance random
quadratic forms. The theory for such quadratic forms is not well
studied. We also refer to an extended version of this paper
\cite{extended} which
covers the class of stochastic volatility processes with regularly
varying marginal distributions.

\section{Preliminaries on the periodogram}\label{sec:1}

The following decomposition of the periodogram is fundamental:
%
\begin{eqnarray}
\label{decomp}
I_{n,X}(\lambda)=\gamma_{n,X}(0) + 2 \sum_{h=1}^{n-1}\cos(\lambda
h) \gamma_{n,X}(h) ,
\end{eqnarray}
where
$
\gamma_{n,X}(h)=\frac1n \sum_{t=1}^{n-|h|}X_tX_{t+h}$, $h\in\mathbb{Z}$,
denotes the {\em sample autocovariance function} of the sample
$X_1,\ldots,X_n$. Note that the definition of $\gamma_{n,X}$
deviates
slightly from the usual one where the $X_t$'s are centered by the sample
mean. However, for the theory given below, this centering is
not essential. Centering with the sample mean $\overline X_n$ is not the
most natural choice when dealing with infinite variance processes.
In what follows, we will frequently make use of the {\em
self-normalized periodogram}
\begin{eqnarray*}
\widetilde I_{n,X}(\lambda)=\dfrac{I_{n,X}(\lambda)}{\gamma
_{n,X}(0)} = \rho_{n,X}(0) + 2 \sum_{h=1}^{n-1}\cos(\lambda h)
\rho_{n,X}(h) ,
\end{eqnarray*}
where
$
\rho_{n,X}(h) = \gamma_{n,X}(h)/\gamma_{n,X}(0)$, $h\in\mathbb{Z}$,
denotes the {\em sample autocorrelation function} of $X_1,\ldots,X_n$.
In view of
(\ref{decomp}), we can rewrite $J_{n,\varepsilon}(f)$ as follows:
%
\begin{equation}
\label{decom2}
J_{n,\varepsilon}(f) = \gamma_{n,\varepsilon}(0)  a_0(f)+ 2 \sum
_{h=1}^{n-1} a_h(f) \gamma_{n,\varepsilon}(h) ,
\end{equation}
where
%
\begin{equation}
\label{eq45}
a_h(f) = \int_0^{\curpi} \cos(\lambda h) f(\lambda)\, \mathrm{d}\lambda
,\qquad h\in\mathbb{Z} ,
\end{equation}
are the {\em Fourier coefficients} of $f$. We also introduce the
self-normalized version of $J_{n,\varepsilon}$:
%
\begin{equation}
\label{decom2a}
\widetilde J_{n,\varepsilon}(f)=\rho_{n,\varepsilon}(0)  a_0(f)+ 2
\sum_{h=1}^{n-1} a_h(f) \rho_{n,\varepsilon}(h) .
\end{equation}
%

\section{The i.i.d. case}\label{sec:iid}

In this section, we study the limit behavior of the integrated
periodograms $J_{n,\varepsilon}$ indexed by classes of functions for an
i.i.d.~S$\alpha$S sequence with
$\alpha\in(0,2)$. In Section~\ref{s1}, we consider the convergence
of the
finite-dimensional distributions. In Sections~\ref{s2} and \ref{s3},
we prove the tightness of
the processes in the cases $\alpha\in(0,1)$ and $\alpha\in[1,2)$,
respectively. In the case $\alpha\in(0,1),$ we solve a more
general weak convergence  problem for random quadratic forms in the
i.i.d.~sequence $(\varepsilon_t)$; the convergence of the integrated
periodograms
indexed by
classes of functions is only a special case. The case $\alpha\in[1,2)$
is more involved. Among others, entropy conditions will be needed and
we only prove results on the weak convergence of the empirical
spectral distribution, that is, we
focus on random quadratic forms with Toeplitz coefficient matrices
given by the Fourier coefficients $a_h(f)$ defined in \eqref{eq45}.

\subsection{Convergence of the finite-dimensional distributions}\label{s1}
A glance at decomposition (\ref{decom2}) is enough to see that
the convergence of the finite-dimensional distributions of
$J_{n,\varepsilon}$ is essentially
determined by the weak limit behavior of the sample
autocovariances $\gamma_{n,\varepsilon}(h)$.
For this reason, we recall a well-known result
due to Davis and Resnick \cite{davisresnick1986}; see also \cite
{brockwelldavis1991}, Section~13.3.
\begin{lemma}\label{le1}
For every $m \ge1$,
%
\begin{equation}
\label{aqua}
\biggl(\dfrac{n \gamma_{n,\varepsilon}(0)}{n^{2/\alpha}} ,\dfrac{n
\gamma_{n,\varepsilon}(h)}{
(n\log n)^{1/\alpha}} ,
h=1,\ldots,m\biggr) \Longrightarrow(Y_0,Y_1,\ldots,Y_m) ,
\end{equation}
where $\Longrightarrow$ denotes weak convergence,
the $Y_h$'s are independent, $Y_0$ is $S_{\alpha/2}(\sigma_1,1,0)$
and $(Y_h)_{h=1,\ldots,m}$ are i.i.d.~$S_{\alpha}(\sigma_2,0,0)$ for
some $\sigma_i=\sigma_i(\alpha)$, $i=1,2$.
In particular,
%
\begin{eqnarray}\label{aqua1}
(n/\log n)^{1/\alpha}
(\rho_{n,\varepsilon}(h))_{h=1,\ldots,m} \Longrightarrow
(Y_h/Y_0)_{i=1,\ldots,m} .
\end{eqnarray}
\end{lemma}

The latter result is an immediate consequence of \eqref{aqua} and
the continuous mapping
theorem.
Lemma~\ref{le1} yields the weak convergence for any finite
linear combination of the sample autocovariances and
autocorrelations. It also suggests that the weak limit
of the standardized process
$J_{n,\varepsilon}(f)$ will be determined by the infinite series
$\sum_{h=1}^{\infty}a_h(f) Y_h$.
But this also means that we need
to require additional assumptions on the sequence $(a_h(f))$.

We will treat this problem in a more general context.
Consider a sequence
\[
\mathbf{a}=(a_1,a_2,\ldots)\in\ell^{\alpha} ,
\]
that is, $\mathbf{a}$ satisfies the
summability condition $\sum_{h}|a_h|^{\alpha}<\infty$.
For such an $\mathbf{a}$, we define the sequences of
processes
%
\begin{equation}\label{eqrto}
\cases{\displaystyle
X_n(\mathbf{a})=(n\log n)^{-1/\alpha}
\sum_{k=1}^{n-1}a_k [n \gamma_{n,\varepsilon}(k)] ,
&\quad$Y(\mathbf{a})= \displaystyle\sum_{k=1}^{\infty} a_k  Y_k ,$
\vspace*{1pt}\cr
\displaystyle\widetilde X_n(\mathbf{a})=
(n/\log n)^{1/\alpha}
\sum_{k=1}^{n-1}a_k \rho_{n,\varepsilon}(k) ,&\quad$\widetilde Y(\mathbf
{a})=Y(\mathbf{a})/Y_0 .$
}
\end{equation}
Here, $Y_0,Y_1,Y_2,\ldots$ are independent stable random variables, as
described
in Lemma~\ref{le1}.
The 3-series theorem \cite{petrov1995}
implies that $\mathbf{a}\in\ell^{\alpha}$ is equivalent to the
a.s.~convergence of the infinite series $Y(\mathbf{a})$
in (\ref{eqrto}). However, for the weak
convergence of $(X_n)$ and $(\widetilde X_n)$, we need a slightly stronger
assumption:
\[
\mathbf{a}\in\ell^{\alpha} \log\ell=\Biggl\{
\mathbf{a}=(a_1,a_2,\ldots)\in\ell^{\alpha}\dvt
\sum_{k=1}^{\infty} |a_k|^{\alpha} \log^{+}\frac{1}{|a_k|}<\infty
\Biggr\} .
\]
This assumption ensures the weak convergence of the random quadratic
forms in (\ref{eqrto}); see the proof of Theorem~\ref{pr1} below.
Assumptions of this type frequently
occur in the literature on infinite variance quadratic forms
(e.g., \cite{kwapienwoyczynski1992}).
They appear in a natural way in tail estimates for quadratic forms in
i.i.d.~stable random variables; see the \hyperref[app]{Appendix}.

We can now formulate our result concerning the convergence of the
finite-dimensional distributions.
\begin{theorem}\label{pr1}
For any
$\alpha\in(0,2)$,
\[
(X_n(\mathbf{a}))_{\mathbf{a}\in\ell^{\alpha} \log\ell}
\stackrel
{\mathit{fidi}}{\longrightarrow} (Y(\mathbf{a}))_{\mathbf{a}\in\ell^{\alpha
} \log\ell}
\quad\mbox{and}\quad(\widetilde X_n(\mathbf{a}))_{\mathbf{a}\in
\ell^{\alpha} \log\ell} \stackrel{\mathit{fidi}}{\longrightarrow}
(\widetilde Y(\mathbf
{a}))_{\mathbf{a}\in\ell^{\alpha} \log\ell} .
\]
\end{theorem}

\begin{pf} Using a Cram\'er--Wold argument, it suffices to
prove the convergence of the one-di\-men\-si\-o\-nal distributions.
From (\ref{aqua}) and the continuous mapping theorem, it immediately follows
that for every $m\ge1$,
\[
(n\log n)^{-1/\alpha}
\sum_{k=1}^{m}a_k [n \gamma_{n,\varepsilon}(k)] \Longrightarrow
Y_m(\mathbf{a})=\sum_{k=1}^{m} a_k  Y_k ,
\]
where $\Longrightarrow$ denotes weak convergence.
Since $\mathbf{a}\in\ell^\alpha$, \mbox{$Y_m(\mathbf{a})
\Longrightarrow Y(\mathbf{a})$} as $m\to\infty$ follows
from the 3-series theorem.
According to \cite{billingsley1968}, Theorem 4.2,
it remains to show that
\begin{eqnarray*}
\lim_{m\to\infty}\limsup_{n\to\infty} P\Biggl((n\log
n)^{-1/\alpha}\Biggl|\sum_{k=m+1}^{n-1} a_k  [n \gamma
_{n,\varepsilon}(k)]\Biggr|>\epsilon\Biggr)=0
\end{eqnarray*}
for every $\epsilon>0$ and $\mathbf{a}\in\ell^{\alpha} \log\ell
$. We write
$p_{n,m}(\mathbf{a};\epsilon)$ for the above probabilities.
Applying Lemma~\ref{le2} in the \hyperref[app]{Appendix}
and the fact that
$\mathbf{a}\in\ell^{\alpha} \log\ell$, we conclude that
\begin{eqnarray*}
p_{n,m}(\mathbf{a};\epsilon)\le \mathit{const}\,  \sum_{k=m+1}^{\infty}
|a_k|^{\alpha} \biggl[1+\log^+ \frac{1}{|a_k|}\biggr]\to0\qquad
\mbox{as $m\to\infty$.}
\end{eqnarray*}
(The constant on the right-hand side depends on $\epsilon$.)
This proves the theorem for $(X_n)$; the convergence
of $(\widetilde X_n)$ can be shown analogously by utilizing (\ref{aqua1}).
\end{pf}

As an immediate corollary of Theorem~\ref{pr1}, we obtain
the following result which solves the problem of finding the
limits
of the finite-dimensional distributions for the integrated periodogram
$J_{n,\varepsilon}$ in (\ref{decom2})
and its self-normalized version $\widetilde J_{n,\varepsilon}$ in
(\ref{decom2a}).
\begin{corollary}\label{co1}
Let $\alpha\in(0,2)$ and
\[
\mathcal{F}=\{f \in L^2[0,\curpi]\dvt  \mathbf
{a}(f)=(a_1(f),a_2(f),\ldots)\in\ell^{\alpha} \log\ell\} ,
\]
where $\mathbf{a}(f)$ is specified in $(\ref{eq45})$.
We then have
\begin{eqnarray*}
n (n \log n)^{-1/\alpha}
[J_{n,\varepsilon}(f)-a_0(f)  \gamma_{n,\varepsilon}(0)
]_{f\in\mathcal{F}}&\stackrel{\mathit{fidi}}{\longrightarrow}&
2 [Y(\mathbf{a}(f))]_{f\in\mathcal{F}} ,
\\
(n/\log n)^{1/\alpha}[\widetilde J_{n,\varepsilon
}(f)-a_0(f)]_{f\in\mathcal{F}}&\stackrel
{\mathit{fidi}}{\longrightarrow}&
2 [\widetilde Y(\mathbf{a}(f))]_{f\in\mathcal{F}} .
\end{eqnarray*}
\end{corollary}

\begin{remark}\label{rem:1}
The condition $\mathbf{a}(f) \in\ell^{\alpha} \log\ell$ is, in
general, not easily
verified. However, if $f$ represents the spectral density of a
stationary process $(X_n)$
with absolutely summable autocovariance function $\gamma_X$,
then, up to a constant multiple, $f$ is represented by the Fourier
series of~$\gamma_X$, and the rate of decay of $\gamma_X(h)\to0$
as $h\to\infty$ is well
known for numerous time series models. For example, if $f$ is the
spectral density of an ARMA process, $\gamma_X(h)\to0$ at an
exponential rate and then $\mathbf{a}(f) \in\ell^{\alpha} \log\ell
$ is satisfied for
every $\alpha>0$.

Conditions ensuring that $\mathbf{a}(f)\in\ell^\alpha$ can be found
in the literature on Fourier series, for example, in \cite
{zygmund2002}. Theorem (3.10) on page 243 in Volume
I of that reference yields, for
Lipschitz continuous functions $f$ with exponent $\beta\in(0,1]$, that
$\mathbf{a}(f)\in\ell^\alpha$ for $\alpha>2/(2\beta+1)$,
but not necessarily for $\alpha=2/(2\beta+1)$. This means, in
particular, that Lipschitz continuous functions do not necessarily satisfy
$\mathbf{a}(f)\in\ell^\alpha$ for small values $\alpha<1$. Zygmund's
Theorem (3.13), \cite{zygmund2002}, page 243, Volume I, states that
$\mathbf{a}(f)\in
\ell^\alpha$
if $f$ is of bounded variation and Lipschitz continuous with exponent
$\beta\in(0,1]$ such that $\alpha>2/(2+\beta)$, but this statement
is not
necessarily valid for $\alpha=2/(2+\beta)$.

We also note that $\mathbf{a}(f)\notin\ell^\alpha$
for $f(x)=I_{[0,x]}$, $x\in(0,\curpi]$ and $\alpha\le1$. Indeed, then
$a_k(f)=k^{-1}\sin(x k)$, $k=1,2,\ldots$ and $\sum_{k}
|a_k(f)|^\alpha=\infty$.
The latter condition implies that the series $Y({\mathbf{a}}(f))$ diverges
a.s.~by the
3-series theorem and the 0--1 law. Hence, Corollary~\ref{co1} does not
apply to the
important class of indicator functions when $\alpha<1$. Moreover,
$(J_{n,\varepsilon}(f))$ is not tight. Indeed, it follows from the argument
above and from \cite{kwapienwoyczynski1992}, Theorem 6.2.1 that for
some $\delta>0$, for every
$K>0$,
\begin{eqnarray*}
\delta&\leq& \lim_{m\to\infty}\lim_{n\to\infty}P\Biggl((n\log
n)^{-1/\alpha}\Biggl|\sum_{k=1}^m a_k(f)
[n \gamma_{n,\varepsilon}(k)]\Biggr|>K\Biggr)
\\
&\le&
\mathit{const}\, \liminf_{n\to\infty} P\bigl(n (n \log n)^{-1/\alpha
}|J_{n,\varepsilon
}(f)-a_0(f) \gamma_{n,\varepsilon}(0)|>K\bigr) .
\end{eqnarray*}
\end{remark}

\subsection{Tightness and weak convergence in the case $\alpha\in
(0,1)$}\label{s2}

In order to derive a full weak convergence counterpart of the
convergence in terms of the finite-dimensional distributions in
Corollary \ref{co1},
it remains to establish tightness of the corresponding family of
laws. We start, once again, in the more general context of random
fields indexed by sequences in $\ell^{\alpha} \log\ell$. Since we
are dealing with the
weak convergence of infinite-dimensional objects, we may expect difficulties
which are due to the geometric properties of the underlying path spaces.
It is also not completely surprising that the case $\alpha\in(0,1)$
is the `better one' in comparison with $\alpha\in[1,2)$;
see, for example, the results on boundedness, continuity and oscillations
of $\alpha$-stable processes
in \cite{samorodnitskytaqqu1994}, Chapter
10. Note, however, that the constraint
$\mathbf{a}(f)\in\ell^{\alpha} \log\ell$ is harder to satisfy
for smaller $\alpha$ than for larger $\alpha$; see also Remark~\ref{rem:1}.

In the present case $\alpha\in(0,1)$, we introduce the function
\[
h(x)=
\cases{
|x|^{\alpha}  \log(b+|x|^{-1}) ,&\quad $x\ne0 ,$
\cr
0 ,&\quad $x=0 ,$
}
\]
where $b$ is chosen so large that $h$ is concave on $(0,\infty)$.
Note that
\[
\ell^{\alpha} \log\ell=\Biggl\{\mathbf{a}\dvt  \sum_{k=1}^{\infty}
h(a_k)  <\infty
\Biggr\}
\]
and this set is a linear metric space when endowed with the
metric
$d(\mathbf{a},\mathbf{b})=\sum_{k=1}^{\infty} h(a_k-b_k).$

Assume that $\mathcal{A}$ is a compact set of $\ell^{\alpha} \log
\ell$ with the
additional property
that
%
\begin{equation}\label{eq2}
\sum_{k=1}^{\infty}
\sup_{\mathbf{a}\in\mathcal{A}} h(a_k)<\infty .
\end{equation}
Observe that $\mathcal{A}$ is then also a compact subset of
$\ell^{\alpha}$ and $(Y(\mathbf{a}))_{\mathbf{a}\in\mathcal{A}}$
is sample-continuous
as a random element with values in $\mathbb{C}(\mathcal{A})$, the
space of
continuous functions on $\mathcal{A}$ equipped with the uniform
topology; see \cite{samorodnitskytaqqu1994},
Section~10.4.
The following is our main result on the weak convergence of
the sequences $(X_n)$ and
$(\widetilde X_n)$
of infinite variance random quadratic forms
in the case $\alpha\in(0,1)$.

\begin{theorem}\label{th1}
Assume $\alpha\in(0,1)$.
For a compact set $\mathcal{A}$ of $\ell^{\alpha} \log\ell$
satisfying $(\ref{eq2}),$
the following weak convergence result holds in $\mathbb{C}(\mathcal{A})$:
\[
(X_n(\mathbf{a}))_{\mathbf{a}\in\mathcal{A}} \Longrightarrow
(Y(\mathbf{a}))_{\mathbf{a}\in\mathcal{A}}\quad\mbox{and}
\quad(\widetilde X_n(\mathbf{a}))_{\mathbf{a}\in\mathcal{A}}
\Longrightarrow(\widetilde Y(\mathbf{a}))_{\mathbf{a}\in\mathcal
{A}} ,
\]
where $X_n$, $\widetilde X_n$, $Y$ and $\widetilde Y$
are defined in $(\ref{eqrto})$ and the processes
$Y$ and $\widetilde Y$ are sample-continuous.
\end{theorem}

\begin{pf} We restrict ourselves to showing that $X_n
\Longrightarrow Y$.
In view of Theorem~\ref{pr1}, it suffices to prove the tightness
of the processes $X_n$ in $(\mathbb{C}(\mathcal{A}),d_{\mathcal
{A}})$, where
$d_{\mathcal{A}}$ is the restriction of $d$ to $\mathcal{A}$.
We have, for positive $\epsilon$ and $\delta$,
%
\begin{eqnarray}\label{eq4a}
&&P\Bigl(\sup_{d_{\mathcal{A}}(\mathbf{a},\mathbf
{b})<\delta} |X_n(\mathbf{a})-X_n(\mathbf{b})|>\epsilon
\Bigr)\nonumber
\\[-8pt]\\[-8pt]
&&\quad \le P\Biggl(\sum_{k=1}^{n-1} \sup_{d_{\mathcal{A}}(\mathbf
{a},\mathbf{b})<\delta}|a_k-b_k| \bigl[n  \gamma_{n,|\varepsilon
|}(k)\bigr]>\epsilon  (n\log n)^{1/\alpha}\Biggr)
= P_n(\epsilon,\delta) .\nonumber
\end{eqnarray}
We want to show that $P_n(\epsilon,\delta)$
can be made arbitrarily small
for all $n$, provided that $\delta$ is small.
We solve this problem in a modified form:
let $\mathbf{C}=(C_0,C_{s,t},s,t=1,2,\ldots)$
be a sequence of i.i.d.~$S_1(1,0,0)$ random variables, independent of
$(\varepsilon_t)$, and $(b_{s,t})$
a double array of real numbers.
We then have
\begin{eqnarray*}
C_0\sum_{1\le s< t\le n} |b_{s,t}|  |\varepsilon_s\varepsilon
_t|\stackrel{d}{=}
\sum_{1\le s< t\le n} b_{s,t} C_{s,t}  |\varepsilon_s\varepsilon_t|
\stackrel{d}{=}\sum_{1\le s< t\le n} b_{s,t} C_{s,t}  \varepsilon
_s\varepsilon_t .
\end{eqnarray*}
By virtue of this argument, it suffices
to replace the products $|\varepsilon_t\varepsilon_s|$
in the quadratic form in (\ref{eq4a}) with the products
$C_{s,t}\varepsilon_t\varepsilon_s$. This
means that it suffices to show that
\[
P_n'(\epsilon,\delta)=
P\Biggl(\sum_{k=1}^{n-1} c_k(\delta) \sum_{j=1}^{n-k} C_{j,j+k}
\varepsilon_j\varepsilon_{j+k}
>\epsilon  (n\log n)^{1/\alpha}\Biggr)
\]
can be made arbitrarily small for all $n$, provided that $\delta$ is
small, where
$c_k(\delta)=\sup_{d_{\mathcal{A}}(\mathbf{a},\mathbf{b})<\delta
}|a_k-b_k|.$
Now, first apply Lemma~\ref{le3} to the $P_n'$'s and then condition
(\ref{eq2}):
\begin{eqnarray*}\label{eq11}
P_n'(\epsilon,\delta)
&\le& \mathit{const}\,  \dfrac{1+\log^+|\epsilon|}{|\epsilon|^{\alpha}}
\frac{1+\log n}{n\log n}   \sum_{i=1}^{n-1}\sum_{k=i+1}^n
|c_k(\delta)|^{\alpha}  \biggl(1+\log^+\frac{1}{|c_k(\delta
)|}\biggr)
\nonumber
\\
&\le&\mathit{const}\,  \sum_{k=1}^{\infty} h(c_k(\delta))\to0\qquad\mbox{as
$\delta\to0$} .
\end{eqnarray*}
\upqed\end{pf}

Theorem~\ref{th1} provides the limit process for a
very general class of random quadratic forms
with infinite first moments. The coefficient matrices of these
quadratic forms are given by Toeplitz matrices.
The conditions on the parameter set $\mathcal{A}$ are nothing but restrictions
on the infinite Toeplitz matrices $(a_{i-j})_{i,j=1,2,\ldots}$.
When specified to the particular case of Fourier coefficients,
as in (\ref{eq45}), Theorem~\ref{th1} yields the following.
\begin{corollary}\label{co2}
Assume that $\alpha\in(0,1)$ and let
\[
\mathcal{F}=\{f\in L^2[0,\curpi]\dvt  \mathbf{a}(f)=(a_1(f),a_2(f),\ldots
)\in\mathcal{A}\} ,
\]
where $\mathcal{A}$ is a compact set of $\ell^{\alpha} \log\ell$
satisfying $(\ref{eq2})$
and $\mathbf{a}(f)$ is specified in $(\ref{eq45})$.
We then have
%
\begin{equation}
\label{refer}
\cases{
n (n \log n)^{-1/\alpha}
[J_{n,\varepsilon}(f)-a_0(f)  \gamma_{n,\varepsilon}(0)
]_{f\in\mathcal{F}}\Longrightarrow
2 [Y(\mathbf{a}(f))]_{f\in\mathcal{F}} ,
\vspace*{3pt}\cr
(n/\log n)^{1/\alpha}[\widetilde J_{n,\varepsilon
}(f)-a_0(f)]_{f\in\mathcal{F}}\Longrightarrow
2 [\widetilde Y(\mathbf{a}(f))]_{f\in\mathcal{F}} ,
}
\end{equation}
where the convergence holds in $\mathbb{C}(\mathcal{F})$.
\end{corollary}

\begin{pf}
Let $T\dvtx  \mathcal{F}\to\mathcal{A}$ be defined by $Tf =
\mathbf{a}(f)$. We claim that $T\mathcal{F}\subset\mathcal{A}$ is
closed, hence
compact. Indeed, if $(f_n)\subset\mathcal{F}$ is such that $Tf_n$
converges in $\ell^{\alpha} \log\ell$ to some point $\mathbf{a}\in
\mathcal{A}$, then
(as $0<\alpha<1$), the sequence
\[
f_n(\lambda) = \frac1\curpi\sum_{j=-\infty}^\infty a_{|j|}(f_n)\cos
j\lambda, \qquad  \lambda\in[0,\curpi] ,
\]
$n=1,2,\ldots$ converges in $L^1[0,\curpi]$ to some function $f$ that
must necessarily be in $\mathcal{F}$. Therefore, $\mathbf{a} = Tf\in
T\mathcal{F}$, and the
latter set is compact.
The above argument shows that $L^2[0,\curpi]$ convergence in
$\mathcal{F}$ is equivalent to $\ell^{\alpha} \log\ell$
convergence in
$T\mathcal{F}$. Since Theorem \ref{th1} implies weak convergence of
the left-hand side of \eqref{refer} to its right-hand side in
$\mathbb{C}(\mathcal{A})$ (when each function $f\in\mathcal{F}$ is
identified with $Tf\in\mathcal{A}$), we conclude that weak
convergence in \eqref{refer} also holds in~$\mathbb{C}(\mathcal{F})$.\looseness=1
\end{pf}\vspace*{1pt}
%

\subsection{Tightness and weak convergence in the case $\alpha\in
[1,2)$}\label{s3}\vspace*{2pt}

Establishing full weak convergence in
the case $\alpha\in[1,2)$ is more difficult than in the case
$\alpha\in(0,1)$.
Indeed, for $\alpha\in(0,1)$, we were allowed to switch from the
random variables $\varepsilon_t$ to their absolute values, due to the
specific geometry
of the spaces $\ell^\alpha$ and, in particular, $\ell^{\alpha} \log
\ell$. The spaces
$\ell^\alpha$, $\alpha\in[1,2)$, have a much more complicated
structure and, therefore, the particular geometry of these spaces will
be need to be invoked in proving tightness for the random quadratic forms
$X_n$ and $\widetilde X_n$. The requirements prescribed by the
geometry are usually
given by entropy conditions;
see \cite{ledoutalagrand1991} for a general
treatment of random elements with values in Banach spaces.
Entropy conditions are typically needed when
$\alpha$-stable processes with $\alpha\in[1,2)$ appear; see
the discussion in \cite{samorodnitskytaqqu1994}, Chapter 12.

In this section, we only consider vectors $\mathbf{a}\in\ell^{\alpha
} \log\ell$
of the form (\ref{eq45}), that is, they are the Fourier coefficients
of some
functions $f$.
Corollary~\ref{co1}
determines the structures of the limit processes
of the quadratic forms
$J_{n,\varepsilon}$ via the convergence of their finite-dimensional
distributions.
It hence suffices to
show the tightness in $\mathbb{C}(\mathcal{F})$ for suitable classes
$\mathcal{F}$.
\cite{kluppelbergmikosch1996a} considered the special case
of the one-dimensional class ${\mathcal{F}}_I$
of indicator functions on $[0,\curpi]$.
We extend their approach to more general classes of
functions, using an entropy condition.

For $f,g\in\mathcal{F}$, let\vspace*{1pt}
\[
d_j(f,g)= j  |a_j(f)-a_j(g)| ,\qquad j\ge1 .
\]
Each $d_j$ defines a pseudo-metric on $\mathcal{F}$.
Let
\[
\rho_k(f,g)=
\max_{2^k\le j<2^{k+1}} d_j(f,g) ,\qquad k\ge0 .
\]
Recall that the $\epsilon$-covering number $N(\epsilon,\mathcal
{F},\rho_k)$
of $(\mathcal{F},\rho_k)$
is the minimal integer $m$ for which we
can find functions $f_1,\ldots,f_m\in\mathcal{F}$ such that
$\sup_{f\in\mathcal{F}} \min_{i=1,\ldots,m}\rho
_k(f,f_i)<\epsilon.$

\begin{theorem}\label{thm3}
Assume that $\alpha\in[1,2)$, define $\mathbf{a}(f)$ as in $(\ref
{eq45})$ and let $\mathcal{F}$ be a subset of $L^2[0,\curpi]$ satisfying:
\begin{longlist}[(ii)]
\item[(i)] $\mathbf{a}(f) \in\ell^{\alpha} \log\ell$ for all $f
\in\mathcal{F}$;
\item[(ii)] $\exists\beta\in(0, \alpha)$ such that
%
\begin{eqnarray}\label{eq:entropy}
N(\epsilon, \mathcal{F}, \rho_k) \le \mathit{const}\,  [1 + (2^k/\epsilon
)^\beta], \qquad
\epsilon> 0,   k \ge0 .
\end{eqnarray}
\end{longlist}
The weak convergence result $(\ref{refer})$ then holds in $\mathbb
{C}(\mathcal{F})$.
\end{theorem}

In contrast to the finite variance case
\cite{dahlhaus1988,mikoschnorvaisa1997},
the entropy condition \eqref{eq:entropy} is a rather strong
one. Indeed, in the papers mentioned, integrability of some power of
$\log N(\epsilon)$
in a neighborhood of the origin suffices. However, conditions
such as \eqref{eq:entropy} are common in problems of continuity and
boundedness for stable processes; see \cite{samorodnitskytaqqu1994},
Chapter 10.

\begin{pf*}{Proof of Theorem \ref{thm3}} The convergence of the finite-dimensional distributions
follows from Theorem~\ref{pr1}.
We restrict ourselves to proving tightness for $J_{n,\varepsilon}$,
which follows
by proving that
%
\begin{eqnarray}\label{eq56}
\lim_{m\to\infty}\limsup_{n\to\infty}
P\Biggl(\Biggl\| \sum_{j=m}^{n} a_j \widehat\gamma_{n,\varepsilon
}(j)\Biggr\|_\mathcal{F}>\epsilon
\Biggr)=0\qquad\mbox{for every $\epsilon>0$,}
\end{eqnarray}
where
$\|g\|_\mathcal{F}=\sup_{f\in\mathcal{F}} |g(f)| $
and
\[
\widehat\gamma_{n,\varepsilon}(j)= (n\log n)^{-1/\alpha} [n \gamma
_{n,\varepsilon}(j)] ,
\qquad j=1,2,\ldots .
\]
As in \cite{kluppelbergmikosch1996a}, (6.4), one can argue that
it suffices in (\ref{eq56}) to consider $m$ and $n$ of some specific
form.
Let $a<b$ be two positive integers,
$m=2^a$ and
$n=2^{b+1}-1$,
and consider numbers
$\epsilon_k=2^{-k\theta}$, $k\ge1,$
with $\theta>0$. For $\theta$ sufficiently small and $a$ large
enough, we
have
\begin{eqnarray}
\label{neu}
P\Biggl( \Biggl\| \sum_{j=m}^{n} a_j \widehat\gamma_{n,\varepsilon}(j)
\Biggr\|_\mathcal{F} >\epsilon\Biggr)
\le  \sum_{k=a}^b P\Biggl( \Biggl\| \sum_{j=2^k}^{2^{k+1}-1} a_j
\widehat\gamma_{n,\varepsilon}(j) \Biggr\|_\mathcal{F} > \epsilon_k
\Biggr)=
\sum_{k=a}^b p_k .
\end{eqnarray}
Consider an array $(\epsilon_{k,l})$ of positive numbers such that
$\epsilon_{k,l}\to0$ as $l\to\infty$ for each $k\ge0$.
Then,
\begin{eqnarray*}
p_k\le
N(\epsilon_{k,0},\mathcal{F},\rho_k)  p_{k,0}
+\sum_{l=1}^{\infty}
N(\epsilon_{k,l},\mathcal{F},\rho_k)  p_{k,l} ,
\end{eqnarray*}

\noindent
where
\begin{eqnarray*}
p_{k,0}&=&\sup_{f\in\mathcal{F}}
P\Biggl(
\Biggl|
\sum_{j=2^k}^{2^{k+1}-1}
a_j(f)  \widehat\gamma_{n,\varepsilon}(j)
\Biggr|> \epsilon_k/2
\Biggr) ,
\\
p_{k,l}&=&
\sup_{f,g\in\mathcal{F} , \rho_k(f,g)\le\epsilon_{k,l-1}}
P\Biggl(\Biggl|
\sum_{j=2^k}^{2^{k+1}-1} [a_j(f)-a_j(g)]
\widehat\gamma_{n,\varepsilon}(j)
\Biggr|
>2^{-(l+1)}\epsilon_k
\Biggr) .
\end{eqnarray*}
By virtue of Lemma~\ref{le2}, we have, for all $f,g\in\mathcal{F}$,
\begin{eqnarray*}
P\Biggl(\Biggl|\sum_{j=2^k}^{2^{k+1}-1}
[a_j(f)-a_j(g)] \widehat\gamma_{n,\varepsilon}(j)
\Biggr|>2^{-(l+1)}\epsilon_k
\Biggr) \le \mathit{const}\, b_{k,l} ,
\end{eqnarray*}
where
\begin{eqnarray*}
b_{k,l}=
\epsilon_k^{-\alpha}2^{\alpha l}
\sum_{j=2^k}^{2^{k+1}-1}|a_j(f)-a_j(g)|^{\alpha}
\bigl[1+ \log^+\bigl(1/|a_j(f)-a_j(g)|\bigr)\bigr] .
\end{eqnarray*}
Assuming that $\rho_k(f,g)\le\epsilon_{k,l-1}$, we have
\begin{eqnarray*}
b_{k,l} &\le& \mathit{const}\, \epsilon_k^{-\alpha}
2^{\alpha l}
\epsilon_{k,l-1}^{\alpha}  \sum_{j=2^k}^{2^{k+1}-1}
j^{-\alpha} [1+\log j\log^+ \epsilon_{k,l-1}^{-1}]
\\
&\le&
\mathit{const}\, \epsilon_k^{-\alpha}
2^{\alpha l}
\epsilon_{k,l-1}^{\alpha} 2^{-k(\alpha-1)}
[1+k \log^+\epsilon_{k,l-1}^{-1}] .
\end{eqnarray*}
Hence, we are left to consider
\begin{eqnarray}
\label{nanuchen}
&&\sum_{k=a}^b \sum_{l=1}^{\infty} N(\epsilon
_{k,l},\mathcal{F},\rho_k)  \epsilon_k^{-\alpha}2^{-k(\alpha-1)
+\alpha l} \epsilon_{k,l-1}^{\alpha} [1+k \log^+\epsilon
_{k,l-1}^{-1} ]\nonumber
\\
&&\quad= \sum_{k=a}^b 2^{-k(\alpha-1-\alpha\theta)} \sum_{l=1}^{\infty}
N(\epsilon_{k,l},\mathcal{F},\rho_k)\epsilon_{k,l-1}^{\alpha}
[1+k \log^+\epsilon_{k,l-1}^{-1}] 2^{\alpha l}
\\
&&\quad\le \mathit{const}\, \sum_{k=a}^b 2^{-k(\alpha-1-\alpha\theta)} \sum
_{l=1}^\infty\biggl[1 + \biggl( \frac{2^k}{\epsilon_{k,l}}
\biggr)^\beta\biggr]\epsilon_{k,l-1}^\alpha[1 + \log^+\epsilon
_{k,l-1}^{-1} ] 2^{\alpha l} .\nonumber
\end{eqnarray}
Assume that $\theta$ is so small that (\ref{neu}) holds. Define the numbers
$
\epsilon_{k,l}=2^{-\gamma_1 l-\gamma_2k}$, $k,l\ge0$
with $\gamma_1,\gamma_2>0$ such that
$1+\gamma_2>(1+\alpha\theta)/(\alpha-\beta)$ and
$\gamma_1>\alpha/(\alpha-\beta)$.
For these parameter choices, it is not difficult to see that
(\ref{nanuchen})
converges to zero by first letting $n\to\infty$ (i.e., $b\to
\infty$) and
then $m\to\infty$ (i.e., $a\to\infty$). This proves (\ref
{eq56}), hence
the tightness of the processes considered in $\mathbb{C}(\mathcal{F})$.
\end{pf*}

In what follows, we give examples of function spaces
$\mathcal{F}$ satisfying condition (ii) of Theorem~\ref{thm3}.
\begin{example}
Consider a space of indexed functions
$\mathcal{G}_\Theta= \{g_\theta\dvt \theta\in\Theta\}$ that
are defined on $[0, \curpi]$ such that $(\Theta, \tau)$ is a
compact metric space, the mapping $\theta\mapsto g_\theta$ is
H\"{o}lder continuous with exponent $b>0$ and constant $K>0$, that is,
\begin{eqnarray*}
\sup_{0 \le x \le\curpi} | g_{\theta_1} (x) - g_{\theta_2} (x)
|
\le K  (\tau(\theta_1, \theta_2))^b \qquad \mbox{for all
} \theta_1, \theta_2 \in\Theta ,
\end{eqnarray*}
and the number of balls (in metric $\tau$) of radius at most $\epsilon
$ necessary to cover $\Theta$ is of the order $\epsilon^{-a}$ for
some $0<a<b\alpha$. Then, $\mathcal{G}_\Theta$ satisfies
$
N(\epsilon, \mathcal{G}_\Theta, \rho_k) \le \mathit{const}\,
(2^k/\epsilon)^{a/b}
$
with $a/b \in(0, \alpha)$. Indeed,
let $\epsilon> 0, k \ge0$. We can find $N \le c ((K\curpi
2^{k+1})/\epsilon)^{a/b}$ balls of radius at most $(\epsilon
/(K\curpi2^{k+1}))^{1/b}$ covering $\Theta$. Call them $B_1, \ldots
, B_N$, with centers $\theta_1, \ldots, \theta_N$. Now, given
$\theta\in\Theta$, we have $\theta\in B_i$ for some $i \in\{
1,\ldots,N\}$ and
\begin{eqnarray*}
\rho_k(g_\theta, g_{\theta_i})
&=& \max_{2^k \le j < 2^{k+1}}j \biggl| \int_0^\curpi\cos(jx) \bigl(
g_\theta(x) - g_{\theta_i}(x) \bigr)\, \mathrm{d}x \biggr|
\\
&\le& 2^{k+1} \curpi\sup_{0 \le x \le\curpi}| g_\theta(x) -
g_{\theta_i}(x) |
\le2^{k+1} \curpi K \tau(\theta, \theta_i)^b
\le\epsilon .
\end{eqnarray*}
The desired bound now follows since
$N(\epsilon, \mathcal{G}_\Theta, \rho_k)
\le N \le \mathit{const}\, (2^k/\epsilon)^{a/b}.$
\end{example}

\begin{example}
Consider a Vapnik--\v Cervonenkis (VC) class $\mathcal{G}$ of functions
defined on $[0,\curpi]$ with VC index $V(\mathcal{G}) = 2$; see
\cite{vandervaartwellner1996}, Section 2.6.2 for more information on
VC classes of
functions. Given $\epsilon> 0$ and $k \ge0$, we can
find $N \le c (\curpi2^{k+1})/\epsilon$ balls of radius at most
$\epsilon/(\curpi2^{k+1})$ that cover $\mathcal{G}$ in the norm
$\frac{1}{\curpi}\int_0^\curpi|\cdot|\, \mathrm{d}x$; see, for example,
\cite{vandervaartwellner1996}, Theorem 2.6.7. Therefore, there exist
$g_1, \ldots, g_N \in\mathcal{G}$ such that for any $g \in\mathcal{G}$,
\begin{eqnarray*}
\min_{1 \le i \le N} \frac{1}{\curpi}\int_0^\curpi|g(x) - g_i(x)|\, \mathrm{d}x <
\frac{\epsilon}{\curpi2^{k+1}} .
\end{eqnarray*}
We then have
\begin{eqnarray*}
\min_{1 \le i \le N} \rho_k(g, g_i) &=& \min_{1 \le i \le N} \max
_{2^k \le j < 2^{k+1}}j \biggl| \int_0^\curpi\cos(jx) \bigl(g(x) - g_i(x)\bigr)
\,\mathrm{d}x \biggr|
\\
&\le& \min_{1 \le i \le N} 2^{k+1} \int_0^\curpi| g(x) - g_i(x)
|\, \mathrm{d}x
\le\epsilon .
\end{eqnarray*}
It follows that
$
N(\epsilon, \mathcal{G}, \rho_k) \le N \le \mathit{const}\,  2^k/\epsilon .
$
\end{example}

\section{The linear process case}\label{sextra}

It is the aim of this section
to show that the results for the case of an i.i.d.~sequence
$(\varepsilon_t)$
translate to the linear process case. The following
decomposition will be crucial:
%
\begin{eqnarray}
\label{munich}
I_{n,X}(\lambda)=I_{n,\varepsilon}(\lambda) |\psi(\mathrm
{e}^{-\mathrm{i}\lambda
})|^2 +R_n(\lambda) .
\end{eqnarray}
This decomposition is analogous to the decomposition $f_{X}(\lambda)=
f_\varepsilon(\lambda) |\psi(\mathrm{e}^{-\mathrm{i}\lambda})|^2$ of
the spectral density $f_X$ of a linear process.
We will show that the normalized
integrated remainder term
$
\int_0^\curpi R_n(\lambda)  f(\lambda)\, \mathrm{d}\lambda
$
is negligible, uniformly over the class of functions $
\mathcal{F}$, in comparison to the normalized main part
\begin{eqnarray*}
\int_0^\curpi
I_{n,\varepsilon}(\lambda)  |\psi(\mathrm{e}^{-\mathrm{i}\lambda})|^2
f(\lambda
) \,\mathrm{d}\lambda ,\qquad f\in\mathcal{F} ,
\end{eqnarray*}
\noindent
which can be treated by the methods of the previous section.
Note that, for a given sequence of coefficients $(\psi_j)_{j\in
{\mathbb{Z}}}$,
the functions $|\psi(\mathrm{e}^{-\mathrm{i}\cdot})|^2  f$ constitute just
another class of functions on $[0,\curpi]$, $\mathcal{F}_{\psi}$ say, and
therefore we will study the process
$
J_{n,\varepsilon}(f)$,
$f\in\mathcal{F_\psi}$,
for suitable classes $\mathcal{F}_\psi$.

\begin{lemma}\label{remainder}
Let $R_n$ be the remainder term appearing in the decomposition \eqref
{munich} of the periodogram $I_{n,X}$. Suppose that the linear filter
$(\psi_j)$ of the process $X$ satisfies
%
\begin{eqnarray}\label{eq35}
\sum_{j=-\infty}^\infty|\psi_j| |j|^{2/\alpha} (1 + \log^+
|j|)^{{(4 - \alpha)/(2\alpha)} + \tau} < \infty
\end{eqnarray}
for some $\tau> 0$ and $\mathcal{F}$ is a collection of real-valued
functions defined on $[0, \curpi]$ such that $\sup_{f \in
\mathcal{F}}\|f\|_2$ $ < \infty$. We then have
\begin{eqnarray*}
\frac{n}{(n \log n)^{1/\alpha}}   \sup_{f \in\mathcal{F}} \biggl|
\int_0^\curpi f(x)R_n(x)\, \mathrm{d}x \biggr| \stackrel{P}{\longrightarrow} 0.
\end{eqnarray*}
\end{lemma}

\begin{pf} From \cite{mikoschgadrich1995},
Proposition 5.1, substituting $n^{1/2}$ for $a_n$, we have the
following decomposition for $R_n$:
%
\begin{equation}
\label{vienna}
R_n(x) = n^{-1} \bigl(\psi(\mathrm{e}^{\mathrm{i}x})L_n(x)K_n(-x) + \psi
(\mathrm{e}
^{-\mathrm{i}x})L_n(-x)K_n(x) + |K_n(x)|^2\bigr) ,
\end{equation}
where $\psi$ is the transfer function as defined before and
\begin{eqnarray*}
L_n(x)&=&\sum_{t=1}^n \varepsilon_t \mathrm{e}^{-\mathrm{i}xt} , \qquad
K_n(x) = \sum_{j=-\infty}^\infty\psi_j \mathrm
{e}^{-\mathrm{i}xj}U_{nj}(x) ,
\\
U_{nj}(x)&=&\Biggl(\sum_{t=1-j}^{n-j} -\sum_{t=1}^n
\Biggr)\varepsilon_t \mathrm{e}^{-\mathrm{i}xt} .
\end{eqnarray*}
We first show that
%
\begin{eqnarray}
\label{izmir}
\frac{1}{(n \log n)^{1/\alpha}}   \sup_{f \in\mathcal{F}} \biggl|
\int_0^\curpi f(x) |K_n(x)|^2\, \mathrm{d}x \biggr| \stackrel{P}{\longrightarrow
} 0 .
\end{eqnarray}
Note that
\begin{eqnarray*}
\biggl| \int_0^\curpi f(x) |K_n(x)|^2 \,\mathrm{d}x \biggr|
&\le& \int_0^\curpi|f(x)| \Biggl(\sum_{j=-\infty}^\infty|\psi_j|
|U_{nj}(x)|\Biggr)^2 \,\mathrm{d}x
\\[-2pt]
&\le&
\mathit{const}\, \Biggl( \sum_{j=-\infty}^{-1} +\sum_{j=1}^\infty\Biggr)
|\psi_j| \int_0^\curpi|f(x)| |U_{nj}(x)|^2 \, \mathrm{d}x .
\end{eqnarray*}
The convergence in (\ref{izmir}) will follow if we can show that the
suprema over $f \in\mathcal{F}$ of the two infinite sums in the last
expression are bounded in probability as $n \rightarrow\infty$. We
will prove this for the second sum; the first one can be handled analogously.

We have, by definition of the terms $U_{nj}(x)$, the Cauchy--Schwarz
inequality
and the fact that, by assumption, $\sup_{f \in\mathcal{F}} \|f\|_2 <
\infty$,
\begin{eqnarray*}
&&\sup_{f\in{\mathcal{F}}}\sum_{j=1}^\infty|\psi_j| \int
_0^\curpi
|f(x)| |U_{nj}(x)|^2 \,\mathrm{d}x
\\[-2pt]
&&\quad\le \sup_{f\in{\mathcal{F}}}\sum_{j=1}^n |\psi_j| \int_0^\curpi|f(x)|
\Biggl|\sum_{t=1-j}^0 \varepsilon_t \mathrm{e}^{-\mathrm{i}xt} - \sum_{t=n-j+1}^n
\varepsilon_t \mathrm{e}^{-\mathrm{i}xt}\Biggr|^2\, \mathrm{d}x
\\[-2pt]
&&\qquad{} + \sup_{f\in{\mathcal{F}}}\sum_{j=n+1}^\infty|\psi_j| \int
_0^\curpi
|f(x)| \Biggl|\sum_{t=1-j}^{n-j} \varepsilon_t \mathrm{e}^{-\mathrm{i}xt} -
\sum
_{t=1}^n \varepsilon_t \mathrm{e}^{-\mathrm{i}xt}\Biggr|^2\, \mathrm{d}x
\\[-2pt]
&&\quad\le c  [I_1(n) + I_2(n) + I_3(n) + I_4(n)] ,
\end{eqnarray*}
where
\begin{eqnarray*}
I_1(n)&=&\sum_{j=1}^n |\psi_j| \Biggl(\int_0^\curpi\Biggl|\sum
_{t=1-j}^0 \varepsilon_t \mathrm{e}^{-\mathrm{i}xt}\Biggr|^4\, \mathrm{d}x \Biggr)^{1/2} ,
\\[-2pt]
I_2(n)&=& \sum_{j=1}^n |\psi_j| \Biggl(\int_0^\curpi\Biggl|\sum
_{t=n-j+1}^n \varepsilon_t \mathrm{e}^{-\mathrm{i}xt}\Biggr|^4\, \mathrm{d}x \Biggr)^{1/2}
,\\[-2pt]
I_3(n)&=& \sum_{j=n+1}^\infty|\psi_j| \Biggl(\int_0^\curpi
\Biggl|\sum_{t=1-j}^{n-j} \varepsilon_t \mathrm{e}^{-\mathrm{i}xt}\Biggr|^4 \,\mathrm{d}x
\Biggr)^{1/2} ,
\\[-2pt]
I_4(n)&=& \sum_{j=n+1}^\infty|\psi_j| \Biggl(\int_0^\curpi\Biggl|
\sum_{t=1}^n \varepsilon_t \mathrm{e}^{-\mathrm{i}xt}\Biggr|^4 \,\mathrm{d}x
\Biggr)^{1/2} .
\end{eqnarray*}
%
It remains to show that each sequence $I_k(n), k=1,2,3,4$, is tight. Now,
\begin{eqnarray*}
I_1(n) \stackrel{d}{=} \sum_{j=1}^n |\psi_j| \Biggl( \int_0^\curpi\Biggl|
\sum_{m=1}^j \varepsilon_m \mathrm{e}^{\mathrm{i}xm} \Biggr|^4 \,\mathrm{d}x
\Biggr)^{1/2} .
\end{eqnarray*}
Let $\epsilon> 0$. Choose $M > 0$ so large that the following holds,
for $\delta= \frac{2\alpha}{4 - \alpha}\tau$:
\begin{eqnarray*}
P\bigl( |\varepsilon_m| > M m^{1/\alpha}(1 + \log m)^{{1/\alpha}
+ \delta} \mbox{ for some } m \ge1 \bigr) \le\epsilon/2 .
\end{eqnarray*}
Write
\begin{eqnarray*}
J_m= \varepsilon_m I_{\{ |\varepsilon_m| \le Mm^{1/\alpha}(1 + \log
m)^{{1/\alpha} + \delta} \}} .
\end{eqnarray*}
Then, for $k > 0$, we have, for $\delta$ chosen as above,
\begin{eqnarray*}
&&P\bigl( I_1(n) > k \bigr)-\epsilon/2
\\
&&\quad\le P \Biggl( \sum_{j=1}^n |\psi_j| \Biggl( \int_0^\curpi
\Biggl| \sum_{m=1}^j \varepsilon_m \mathrm{e}^{\mathrm{i}tx} \Biggr|^4  \,\mathrm{d}x
\Biggr)^{1/2} \Biggl( \int_0^\curpi\Biggl| \sum_{m=1}^j J_m  \mathrm
{e}^{\mathrm{i}tx}
\Biggr|^4 \,\mathrm{d}x \Biggr)^{1/2} > k \Biggr)
\\
&&\quad \le k^{-1} \sum_{j=1}^n |\psi_j|
\Biggl( \int_0^\curpi E \Biggl| \sum_{m=1}^j J_m
\mathrm{e}^{\mathrm{i}tx} \Biggr|^4 \, \mathrm{d}x \Biggr)^{1/2}
\\
&&\quad = k^{-1} \sum_{j=1}^n |\psi_j|
\Biggl[\int_0^\curpi\Biggl(\sum_{m=1}^j E (J_m^4)
\\
&&{}\qquad + 6  \sum_{m_1=1}^j \sum_{m_2=m_1+1}^j E (
J_{m_1}^2 ) E(J_{m_2}^2) \cos\bigl((m_1-m_2)x\bigr) \,\mathrm{d}x\Biggr) \Biggr]^{1/2}
\\
&&\quad \le \frac{c}{k} \sum_{j=1}^n |\psi_j|\Biggl[ \Biggl(\sum_{m=1}^j
\bigl(
m^{1/\alpha}(1 + \log m)^{{1/\alpha} + \delta} \bigr)^{4 -
\alpha} \Biggr)^{1/2}
\\
&&{}\qquad + \Biggl( \Biggl[\sum_{m=1}^j \bigl( m^{1/\alpha
}(1 + \log m)^{{1/\alpha} + \delta} \bigr)^{2 - \alpha}
\Biggr]^2 \Biggr)^{1/2}\Biggr]
\\
&&\quad \le \frac{c}{k} \sum_{j=1}^\infty|\psi_j|  j^{2/\alpha}  (1 +
\log j)^{{(4-\alpha)/(2\alpha)} + \tau} .
\end{eqnarray*}
By virtue of \eqref{eq35}, the last expression can be made smaller
than $\epsilon/2$ by choosing $k$ large enough, which proves the
tightness of $I_1(n)$. Similar arguments show that $I_j(n)$,
$j=2,3,4,$ are tight sequences as well. The convergence in (\ref
{izmir}) follows.

By the decomposition (\ref{vienna}), the proof will be finished if we
can also establish that
%
\begin{equation}
\label{prague}
\frac{1}{(n \log n)^{1/\alpha}}   \sup_{f \in\mathcal{F}} \biggl|
\int_0^\curpi f(x)\psi(\mathrm{e}^{\mathrm{i}x})L_n(x)K_n(-x)\, \mathrm{d}x \biggr|
\stackrel
{P}{\longrightarrow} 0 .
\end{equation}
We have, by the Cauchy--Schwarz inequality and the identity $|L_n(x)|^2
= n I_{n,\varepsilon}(x)$,
\begin{eqnarray*}
&&\biggl| \int_0^\curpi f(x)\psi(\mathrm{e}^{\mathrm{i}x})L_n(x)K_n(-x)\, \mathrm{d}x\biggr|
\\
&&\quad \le c \|f\|_2 \biggl( \int_0^\curpi|L_n(x)K_n(-x)|^2 \,\mathrm{d}x
\biggr)^{1/2}\\
&&\quad \le c \|f\|_2   n^{1/2} \Bigl( \sup_{0 \le x \le\curpi
}I_{n,\varepsilon}(x) \Bigr)^{1/2} \biggl(\int_0^\curpi
|K_n(-x)|^2\, \mathrm{d}x \biggr)^{1/2} .
\end{eqnarray*}
We therefore see that
\begin{eqnarray*}
&&\frac{1}{(n \log n)^{1/\alpha}}   \sup_{f \in\mathcal
{F}} \biggl| \int_0^\curpi f(x)\psi(\mathrm{e}^{\mathrm{i}x})L_n(x)K_n(-x)\, \mathrm{d}x
\biggr|
\\
&&\quad \le \frac{c}{n^{{1/\alpha} - {1/2}}}   \frac{
( \sup_{0 \le x \le\curpi}I_{n,\varepsilon}(x) )^{1/2}}{(\log
n)^{1/\alpha}} \biggl(\int_0^\curpi|K_n(-x)|^2 \,\mathrm{d}x
\biggr)^{1/2} .
\end{eqnarray*}
Similar arguments as for (\ref{izmir}) ensure the tightness of
the sequence $\int_0^\curpi|K_n(-x)|^2 \,\mathrm{d}x$. The tightness
of the term
\begin{eqnarray*}
\frac{( \sup_{0 \le x \le\curpi}I_{n,\varepsilon}(x)
)^{1/2}}{(\log n)^{1/\alpha}}
\end{eqnarray*}
follows from \cite{mikoschresnic2000} Theorem 2.1
(for $0 < \alpha< 1$) and
Proposition 3.1 (for $1 \le\alpha< 2$).
Thus, we conclude that (\ref{prague}) holds,
and Lemma \ref{remainder} is proved.
\end{pf}

\begin{remark}
A referee kindly pointed out that Lemma~\ref{remainder} remains valid
under the following condition, which is weaker than \eqref{eq35}:
assume that there exists a sequence $(\omega_n)$ of positive numbers
such that
\
\begin{eqnarray*}
(n\log
n)^{-1/\alpha} \Biggl(\sum_{k=1}^n\omega_k^{-\alpha}
\Biggr)^{2/\alpha-1/2}
\sum_{|j|\le n} |\psi_j| \Biggl(\sum_{l=1}^j\omega_l^{4-\alpha
}\Biggr)^{1/2}\to0 .
\end{eqnarray*}
Condition \eqref{eq35} follows by taking $\omega_n=n^{1/\alpha}
(1+ \log n)^{1/\alpha+\delta}$ for some positive $\delta$.
\end{remark}

By (\ref{munich}), we may write, for each $f$,
\begin{eqnarray*}
J_{n,X}(f) - a_0(f|\psi|^2)\gamma_{n,\varepsilon}(0) &=&
J_{n,\varepsilon}(f|\psi|^2) - a_0(f|\psi|^2)\gamma_{n,\varepsilon
}(0)
\\
&&{}+ \int_0^\curpi f(x) R_n(x)\, \mathrm{d}x ,
\end{eqnarray*}
where $|\psi|^2$ stands for $|\psi(\mathrm{e}^{-\mathrm{i}\cdot})|^2$. Combining this
decomposition with Lemma \ref{remainder}, we can now state the
following analogs to Corollary \ref{co2} and Theorem \ref{thm3}.\vspace*{1pt}

\begin{corollary}
Assume that $\alpha\in(0,1)$ or $\alpha\in[1,2)$ and let $\mathcal
{F}$ be as defined
as in Corollary \ref{co2} or Theorem~\ref{thm3}, respectively.
Suppose that the set $\{f\dvtx [0, \curpi] \rightarrow\mathbb{R}\dvt f|\psi|^2
\in\mathcal{F}\} = \mathcal{F}_\psi$ satisfies
$\sup_{f \in\mathcal{F}_\psi}\|f\|_2 < \infty$ and \eqref{eq35}
holds for some $\tau> 0$.
We then have\vspace*{2pt}
%
\begin{equation}\label{copenhagen}
\cases{
n (n \log n)^{-1/\alpha}
[J_{n,X}(f)-a_0(f|\psi|^2)  \gamma_{n,\varepsilon}(0)
]_{f\in\mathcal{F}_\psi}\Longrightarrow
2 [Y(\mathbf{a}(f|\psi|^2))]_{f\in\mathcal{F}_\psi} ,
\vspace*{3pt}\cr
(n/\log n)^{1/\alpha}[\widetilde J_{n,X}(f)-a_0(f|\psi
|^2)]_{f\in\mathcal{F}_\psi}\Longrightarrow
2 [\widetilde Y(\mathbf{a}(f|\psi|^2))]_{f\in\mathcal{F}_\psi} ,
}\vspace*{2pt}
\end{equation}
where the convergence holds in $\mathbb{C}(\mathcal{F_\psi})$.
\end{corollary}\vspace*{1pt}

\begin{appendix}
\section*{Appendix}\label{app}\vspace*{1pt}

For an array $\mathbf{b}=(b_{s,t})$ of real numbers, define the quadratic
forms\vspace*{2pt}
\[
Q_{n,\varepsilon}(\mathbf{b})=\sum_{1\le s\ne t\le n}b_{s,t}
\varepsilon_s\varepsilon_t\vspace*{2pt}
\]
and\vspace*{2pt}
\[
\Gamma_n(\mathbf{b})=
\sum_{1\le s\ne t\le n} |b_{s,t}|^{\alpha} \biggl(1+\log^+ \frac{1}
{|b_{s,t}|}\biggr) .\vspace*{2pt}
\]
The following lemma is a consequence of \cite
{rosinskiwoyczynski1987}, Theorem~3.1; see also
\cite{kwapienwoyczynski1992}.\vspace*{1pt}

\renewcommand{\thelemma}{A.\arabic{lemma}}
\setcounter{lemma}{0}
\begin{lemma}\label{le2}
For $\alpha\in(0,2)$, there exists a positive constant $D_{\alpha}$
such that for all $x>0$,\vspace*{2pt}
\begin{eqnarray*}
P\bigl(Q_{n,\varepsilon}(\mathbf{b})>x\bigr)\le
D_{\alpha} \dfrac{1+\log^+x}{ x^{\alpha}} \Gamma_n(\mathbf{b}) .\vspace*{2pt}
\end{eqnarray*}
\end{lemma}

Now, let $\mathbf{C}=(C_0,C_{s,t},s,t=1,2,\ldots)$
be a sequence of i.i.d.~$S_1(1,0,0)$ random variables, independent of
$(\varepsilon_t)$, and let $\mathbf{b}$
be as above.
The following lemma is a consequence of Lemma~\ref{le2}.\looseness=1

\begin{lemma}\label{le3}
For $\alpha\in(0,1)$, there exists a positive
constant $D_{\alpha}'$ such that for all $x>0$,
\renewcommand{\theequation}{A.\arabic{equation}}
\setcounter{equation}{0}
\begin{equation}
\label{eq18}
I_n(x) = P\biggl(\sum_{1\le s< t\le n} b_{s,t} C_{s,t}  \varepsilon
_s\varepsilon_t>x\biggr) \le D_{\alpha}' \dfrac{1+\log
^+x}{x^{\alpha}} \Gamma_n(\mathbf{b}) .
\end{equation}
\end{lemma}

\begin{pf}
Apply Lemma~\ref{le2} to $I_n(x)$, conditionally on
$\mathbf{C} $:
\begin{eqnarray}
\label{eq8}
I_n(x)
&=& P\biggl(\sum_{1\le s< t\le n} b_{s,t} C_{s,t}  \varepsilon
_s\varepsilon_t>x\biggr)\nonumber
\\
&=& E_\mathbf{C} P\biggl(\sum_{1\le s< t\le n} b_{s,t}
C_{s,t}  \varepsilon_s\varepsilon_t>x \Big|\mathbf{C}\biggr)
\\
&\le& \mathit{const}\,  \dfrac{1+\log^+ x}{x^{\alpha}}   \sum_{s=1}^n\sum
_{t=i+1}^n |b_{s,t}|^{\alpha}   E|C_{0}|^{\alpha}  \biggl(1+\log^+
\frac{1}{|b_{s,t} C_{0}|}\biggr), \nonumber
\end{eqnarray}
Because $\alpha\in(0,1)$, we also have, for $x>0$,
%
\begin{equation}
\label{eq9}
E|C_{0}|^{\alpha} \biggl(1+\log^+ \frac{1}{|x C_{0}|}\biggr)
\le \mathit{const}\,   \biggl(1+\log^+ \frac{1}{|x|}\biggr) ,
\end{equation}
and combining (\ref{eq8}) and (\ref{eq9}), we thus obtain (\ref{eq18}).
\end{pf}

\end{appendix}

\section*{Acknowledgements}
Thomas Mikosch's research was
supported in part by the Danish Research Council (FNU) Grant
272-06-0442. Gennady Samorodnitsky's research was supported in part by
the ARO Grant W911NF-07-1-0078 at Cornell University and a Villum Kann
Rasmussen Visiting Professor
Grant at the University of Copenhagen. The paper was written when Sami
Umut Can
and Gennady Samorodnitsky visited the Department of Mathematical
Sciences at the University of Copenhagen in 2008/2009. They take pleasure
in thanking
the Department for hospitality and financial support. The authors
also thank the anonymous referees for suggestions that helped to
improve the paper.

\printhistory

\end{document}